For any comments or questions, please contact ma_mike@tongji.edu.cn

# A Critical Note on "Information Distortion in a Supply Chain: The Bullwhip Effect"

## *On Order Batching*

Hau Mike Ma, Jiazhen Huo

**Abstract.**

In the seminal paper "Information Distortion in a Supply Chain: The Bullwhip Effect" (Lee, et al. 1997, hereafter referred to as LPW), order batching is regarded as one of the four sources of the bullwhip effect. LPW proved that, in all cases (random ordering, balanced ordering, and correlated ordering), order batching will surely lead to bullwhip effects. However, we identify two improper assumptions in LPW. First, the batched order $Z_t$ is de facto the moving summation of previous demands, including overlapping demands. In fact, the batched order should be modeled as periodic summation of previous demands. Second, in the random ordering case, the number of retailers *n* is modeled as a binomial variable which is identically distributed for a randomly chosen period *t* in a review cycle. In fact, *n* should follow a sequential hypergeometric distribution. To address the two issues, we decompose a demand sequence using law of the total variance, exploring variance interplay between batched and non-batched demands in the positively correlated case. We find that even under the contrived i.i.d. assumption in LPW, order batching does not necessarily lead to the bullwhip effect.

1. **Introduction**

Order batching could lead to significant bullwhip effect based on Formulas (3.11-3.13) in LPW (Lee et al. 1997). Take positively correlated ordering case as example. Suppose there are 2 retailers ($N = 2$), each using a periodic review system with the review cycle of 2 periods ($R = 2$). Assume that mean and variance of demands for each retailer are 10 and 1 ($m = 10, \sigma^2 = 1$). According to LPW Formula (3.12), the variance after order batching is $N\sigma^2 + m^2 N^2 (R - 1)$, which is 402. In contrast, the variance of non-batched demands $N\sigma^2$ is 2, leading to a substantial bullwhip ratio of 201. This magnitude of bullwhip ratio is inconsistent with empirical findings, which typically report single-digit values (Yao et al. 2020). This raises the question: why does the bullwhip effect due to order batching appear significantly larger in LPW's analysis compared to empirical observations?

In Section 2, we identify two inappropriate assumptions in the mathematical derivation of LPW that can cause this discrepancy. In Section 3, we present refined formulas to calculate the variance of batched orders.

2. **Ergodicity of the batched order $Z_t$ and statistic distribution of $n$**

Consider there exists one retailer. Three ordering cases in LPW are degraded to the same scenario. Now define the batched order of one retailer as $D_t$ and $D_t := \sum_{k=t-R}^{t-1} \xi_k$. If $D_t$ is ergodic (Hamilton 1994, p.46), it means that $D_t$ is the moving summation of previous $R$ demands ($R \geq 2$), which means it contains overlapping demands. For example, assume $R = 2$, we can get $D_3 = \xi_2 + \xi_1$ and $D_4 = \xi_3 + \xi_2$, thus $\xi_2$ is double counted in the batched orders. This could lead to biased impact on the variance of batched orders. However, based on the supply chain mechanism in LPW, $D_t$ should be the periodic summation of $R$ demands of the previous review cycle. At each review period where $t \bmod R = 1$, $D_t$ is the sum of demands of the previous $R$ periods. For periods where $t \bmod R \neq 1$, $D_t = \emptyset$. Therefore, if $D_t$ is ergodic at time scale of $t$, the batched order $D_t$ will include overlapping demands, contradictory to the ordering system designed in LPW. If $D_t$ is not ergodic at time scale of $t$, it's incorrect to treat $D_t$ as a random variable representing system properties of batched orders. In summary, mathematical methods used to construct periodic demand batching in Formulas (3.11-3.13) of LPW are not appropriate.

Now consider there are $N$ retailers. In the random ordering case presented in LPW, $n$ is defined as a random variable denoting the number of retailers placing orders in a randomly chosen period. According to LPW, $n$ is modeled as a binomial variable that follows a binomial distribution, implying that $n_t$ (with subscript $t$ newly introduced) is identical within a review cycle. There arises the possibility that the same retailer could place orders multiple times in one review cycle. However, in the random ordering case, $n_t$ should follow a sequential hypergeometric distribution. Consequently, $Z_{t_1}, \ldots, Z_{t_R}$ in one review period are also different in their statistic properties.

To summarize, in random ordering case, $Z_t$ is not identical for all time periods. In other words, ergodicity of $Z_t$ for expectation and variance is not guaranteed at time scale of $t$. In the positively correlated ordering case, $n$ is modeled to be a two-outcome random variable, which also makes $Z_t$ unidentical within one review cycle. This issue also exists in the balanced ordering case.

As ergodicity of the batched order $Z_t$ can't be guaranteed at time scale $t$, $Z_t$ can't be treated as a random variable with ordering system properties. Therefore, Formula (3.11-3.13) in LPW should be reexamined. In the next section, we model the demand system as a sequence to analyze the effects of order batching. Our findings indicate that order batching does not necessarily result in the bullwhip effect.

## 3. Rewriting formulas in the positively correlated ordering case

In this section, we focus on the positively correlated case, as the batched orders in this scenario can be clearly divided into two steps, each governed by distinct mechanisms. First, we calculate the variance of batched orders resulting from time aggregation for a single retailer. Second, we extend the analysis to calculate the variance of batched orders across $N$ retailers. This approach enables the analysis of batched orders from two perspectives: order batching by period and order batching by retailer.

### 3.1 Calculating variance of batched demands under time aggregation for one retailer

Consider a single retailer with a demand sequence of size $T$ from the random variable $\Xi_t$, denoted as $\{\xi_t\}_{t=1}^T = \{\xi_1, \xi_2, \ldots, \xi_T\}$. Divide this sequence by $R$ into $M$ review cycles, such that $T = R \cdot M$. For simplicity, assume both $M$ and $R$ are integers. Let the demands in one review cycle be represented as $\{\xi_t^{(i)}\}_{t=1}^R$ for $i = 1, 2, \ldots, M$. The mean of the $i$-th review cycle is denoted as $\bar{\xi}^{(i)} = \frac{1}{R}\sum_{t=1}^R \xi_t^{(i)}$ and the variance of the $i$-th review cycle is denoted as $\sigma_{\xi_t^{(i)}}^2 = \frac{1}{R}\sum_{t=1}^R \left(\xi_t^{(i)} - \bar{\xi}^{(i)}\right)^2$. Define $\sigma_{\xi_{within}}^2$ as the expectation of demand variances of all review cycles, and $\sigma_{\xi_{within}}^2 = \frac{1}{M}\sum_{i=1}^M \sigma_{\xi_t^{(i)}}^2$. Define $\sigma_{\bar{\xi}^{(i)}}^2$ as the variance of means from all review cycles. The variance of demands of the retailer is denoted as $\sigma_{\xi_{total}}^2$.

According to the law of total variance, the total variance can be expressed as the sum of the expected variance of demand within review cycles and the variance of the mean demand across review cycles.

Thus,

$$\sigma_{\xi_{total}}^2 = \sigma_{\xi_{within}}^2 + \sigma_{\bar{\xi}^{(i)}}^2 \tag{1}$$

We denote the batched demands of the $i$-th review cycle as $\xi_{agg}^{(i)}$, and $\xi_{agg}^{(i)} = \sum_{t=1}^R \xi_t^{(i)}$, which is the demand faced by the supplier. Let $\bar{\xi}^{(i)} = \frac{1}{R}\sum_{t=1}^R \xi_t^{(i)}$ be the mean of the $i$-th review cycle of demands. We can get $\xi_{agg}^{(i)} = R\bar{\xi}^{(i)}$. Let $\sigma_{\xi_{agg}^{(i)}}^2$ denote the variance of the batched demands across $M$ review cycles and $\sigma_{\bar{\xi}^{(i)}}^2$ denote the variance of averaged orders across $M$ review cycles, the relationship between the variance of batched demands and the variance of averaged demands is

$$\sigma_{\xi_{agg}^{(i)}}^2 = R^2 \sigma_{\bar{\xi}^{(i)}}^2 \tag{2}$$

Combining Formula (1) and (2), we can get

$$\sigma_{\xi_{agg}^{(i)}}^2 = R^2(\sigma_{\xi_{total}}^2 - \sigma_{\xi_{within}}^2) \tag{3}$$

As $T \to \infty$, $\sigma_{\xi_{total}}^2$ converges to $\sigma^2$, which means demand variance of each retailer is $\sigma^2$ over time. Thus, formula (3) can be written as

$$\sigma_{\xi_{agg}^{(i)}}^2 = R^2(\sigma^2 - \sigma_{\xi_{within}}^2) \tag{3.1}$$

### 3.2 Calculating the variance of batched order for $N$ retailers

Suppose there exist $N$ retailers each using a periodic review system with the review cycle equal to $R$ periods. We add a subscript $j$ for retailers. Formula (3) can be written as

$$\sigma_{\xi_{agg,j}^{(i)}}^2 = R^2\left(\sigma_{\xi_{total,j}}^2 - \sigma_{\xi_{within,j}}^2\right) \tag{3.2}$$

Let $Z_i$ denote the batched order of $N$ retailers for the $i$-th review cycle and $Z_i = \sum_{j=1}^N \xi_{agg,j}^{(i)}$. $Z_i$ is ergodic at time scale of a review cycle, and $\xi_{agg,j}^{(i)}$ are also ergodic at time scale of a review cycle. Therefore, $Var(Z_i) = \sum_{j=1}^N \sigma_{\xi_{agg,j}^{(i)}}^2$. Hence, the variance of batched orders across review cycles for $N$

retailers is

$$Var(Z_i) = \sum_{j=1}^{N} \sigma^2_{\xi_{agg,j}^{(i)}} = \sum_{j=1}^{N} \left[ R^2 \left( \sigma^2_{\xi_{total,j}} - \sigma^2_{\xi_{within,j}} \right) \right] \tag{4}$$

As $T \to \infty$, $\sigma^2_{\xi_{total,j}}$ converges to $\sigma^2$, which means demand variance of each retailer is $\sigma^2$ over time. Thus, formula (4) can be written as

$$Var(Z_i) = \sum_{j=1}^{N} \left[ R^2 \left( \sigma^2 - \sigma^2_{\xi_{within,j}} \right) \right] = R^2 N \sigma^2 - R^2 \sum_{j=1}^{N} \sigma^2_{\xi_{within,j}} \tag{4.1}$$

### 3.3 Compare the variance of batched orders and the variance of non-batched demands

#### 3.3.1 Order batching for one retailer

When $N = 1$, we compare Formula (3.1) with $\sigma^2$, the difference is

$$R^2(\sigma^2 - \sigma^2_{\xi_{within}}) - \sigma^2$$

Which can be written as

$$\sigma^2(R^2 - 1) - \sigma^2_{\xi_{within}} R^2$$

Denominate $\sigma^2(R^2 - 1) - \sigma^2_{\xi_{within}} R^2$ with $\sigma^2 R^2$, we get

$$\frac{R^2 - 1}{R^2} - \frac{\sigma^2_{\xi_{within}}}{\sigma^2} \tag{5}$$

Therefore, compare $\frac{R^2-1}{R^2} - \frac{\sigma^2_{\xi_{within}}}{\sigma^2}$ with 0, we will know whether order batching with period aggregation will cause bullwhip effect. Here are three scenarios:

**Scenario A:** When $\frac{\sigma^2_{\xi_{within}}}{\sigma^2} = \frac{R^2-1}{R^2}$, variance of batched orders is the same to variance of demands experienced by the retailer, meaning order batching have no effect on demand variance.

**Scenario B:** When $\frac{\sigma^2_{\xi_{within}}}{\sigma^2} > \frac{R^2-1}{R^2}$, variance of batched orders is less than the variance of demands experienced by the retailer, meaning order batching dampens demand variance.

**Scenario C:** When $\frac{\sigma^2_{\xi_{within}}}{\sigma^2} < \frac{R^2-1}{R^2}$, variance of batched orders is larger than the variance of demands experienced by the retailer, meaning order batching amplifies demand variance, causing bullwhip effect.

#### 3.3.2 Order batching for $N$ retailer

If $N > 1$, we compare Formula (4.1) with $N\sigma^2$, the difference is

$$R^2 N \sigma^2 - R^2 \sum_{j=1}^{N} \sigma^2_{\xi_{within,j}} - N\sigma^2$$

Which can be written as

$$(R^2 - 1)N\sigma^2 - R^2 \sum_{j=1}^{N} \sigma^2_{\xi_{within,j}}$$

Denominate $(R^2 - 1)N\sigma^2 - R^2 \sum_{j=1}^{N} \sigma^2_{\xi_{within,j}}$ with $N\sigma^2 R^2$, we get

$$\frac{R^2-1}{R^2} - \frac{\sum_{j=1}^{N}\sigma^2_{\xi_{within,j}}}{N\sigma^2} \qquad (6)$$

Similarly, Scenarios A through C also hold, leading to the conclusion that positively correlated ordering does not necessarily result in the bullwhip effect.

## 4. Conclusion

Order batching, identified as a primary driver of the bullwhip effect in LPW, warrants reevaluation due to two questionable assumptions. First, the variance of batched orders, when treated as moving sum of the demands from the previous review cycle, incorporates overlapping demands, which contradicts the periodic batching mechanism. Second, the random variable $n$ is not identical across different time periods. To address these issues in LPW, we model the demand system as a demand sequence and apply the law of total variance to examine the effects of order batching in positively correlated ordering case. Our findings demonstrate that order batching does not necessarily lead to the bullwhip effect, even under contrived assumptions in LPW. Given the prevalence practice of order batching in industries, our results shed light on the potential of achieving economies of scale while mitigating the bullwhip effect.

For any comments or questions, please contact ma_mike@tongji.edu.cn